\documentclass[12pt]{article}

\makeatletter
\renewcommand{\@oddfoot}{\hfill \thepage}
\makeatother

\usepackage[T2A]{fontenc}
\usepackage[cp1251]{inputenc}
\usepackage[russian, english]{babel}
\usepackage{mathrsfs}
\usepackage{amssymb}
\usepackage{amsmath,amsthm,amsfonts}
\usepackage{cite}
\usepackage[dvipdf]{graphicx}
\usepackage{fancyhdr}

\begin{document}

\begin{center}
{\bf INTEGRAL EQUATION FOR THE TRANSITION DENSITY \\ OF THE MULTIDIMENSIONAL MARKOV RANDOM FLIGHT}
\end{center}

\begin{center}
Alexander D. KOLESNIK\\
Institute of Mathematics and Computer Science\\
Academy Street 5, Kishinev 2028, Moldova\\
E-Mail: kolesnik@math.md
\end{center}

\vskip 0.2cm

\begin{abstract}
We consider the Markov random flight $\bold X(t)$ in the Euclidean space $\Bbb R^m, \; m\ge 2,$ 
starting from the origin $\bold 0\in\Bbb R^m$ that, 
at Poisson-paced times, changes its direction at random according to arbitrary distribution on the unit 
$(m-1)$-dimensional sphere $S^m(\bold 0,1)$ having absolutely continuous density. 
For any time instant $t>0$, the convolution-type recurrent relations for the joint and conditional densities of 
process $\bold X(t)$  and of the number of changes of direction, are obtained. Using these relations, we derive 
an integral equation for the transition density of $\bold X(t)$ whose solution is given in the form of a uniformly 
converging series composed of the multiple double convolutions of the singular component of the density with itself. 
Two important particular cases of the uniform distribution on $S^m(\bold 0,1)$ and of the Gaussian distributions on the 
unit circumference $S^2(\bold 0,1)$ are separately considered. 
\end{abstract}

\vskip 0.1cm


\vskip 0.2cm


\section{Introduction}

\numberwithin{equation}{section}

Random motions at finite speed in the multidimensional Euclidean spaces $\Bbb R^m, \; m\ge 2$, also called random flights, have become the subject of intense researches in last decades. The majority of published works deal with the case of isotropic Markov random flights when the motions are controlled by a homogeneous Poisson process and their directions are taken uniformly on the unit $(m-1)$-dimensional sphere \cite{kol2,kol3,kol4,kol5,kol6,kol7}, \cite{mas}, \cite{sta1,sta2}. The limiting behaviour of the Markov random flight with a finite number of fixed directions in $\Bbb R^m$ was examined in \cite{ghosh}. 
In recent years the non-Markovian multidimensional random walks with Erlang- and Dirichlet-distributed displacements were studied in a series of works \cite{lecaer1,lecaer2,lecaer3,letac}, \cite{pogor1,pogor2}. Such random motions at finite velocities are of a great interest due to both their big theoretical importance and numerous fruitful applications in physics, chemistry, biology and other fields. 

When studying such a motion, its explicit distribution is, undoubtedly, the most attractive aim of the researches. However, despite many efforts, the closed-form expressions for the distributions of Markov random flights were obtained in a few cases only. In the spaces of low even-order dimensions such distributions were obtained in explicit forms by different methods (see \cite{sta2}, \cite{mas}, \cite{kol5}, \cite{kol7} for the Euclidean plane $\Bbb R^2$, \cite{kol6} for the space $\Bbb R^4$ and \cite{kol2} for the space $\Bbb R^6$). Moreover, in the spaces $\Bbb R^2$ and $\Bbb R^4$ such distributions are surprisingly expressed in terms of elementary functions, while in the space $\Bbb R^6$ the distribution has the form of a series composed of some polynomials. As far as the random flights in the odd-dimensional Euclidean spaces is concerned, their analysis is much more complicated in comparison with the even-dimensional cases. A formula for the transition density of the symmetric Markov random flight with unit speed in the space $\Bbb R^3$ was given in \cite{sta1}, however it has a very complicated form of an integral with variable limits whose integrand depends on the inverse tangent function from the integration variable (see \cite[formulas (1.3) and (4.21)]{sta1}). Moreover, the density presented in this work evokes some questions since its absolutely continuous (integral) part is discontinuous at the origin $\bold 0\in\Bbb R^3$ and this fact seems fairly strange. 

The characteristic functions of the multidimensional random flights are much more convenient objects for analysing than their densities. 
This is due to the fact that, while the densities are finitary fucntions (that is, functions defined on the compact sets of $\Bbb R^m$), their characteristic functions (Fourier transforms) are analytical real functions defined everywhere in $\Bbb R^m$. That is why just the characteristic functions have become the subject of the vast research whose results were published in \cite{kol3} and \cite{kol8}. In particular, in \cite{kol3} the tme-convolutional recurrent relations for the joint and conditional characteristic functions of the Markov random flight in the Euclidean space $\Bbb R^m$ of arbitrary dimension $m\ge 2$, were obtained. By using these recurrent ralations, the Volterra integral equation of second kind with continuous kernel for the unconditional characteristic function was derived and a closed-form expression for its Laplace transform was given. Such convolutional structure of the characteristic functions makes hints about a similar one  for the respective densities. Discovering of such convolutional relations for the densities of Markov random flights in $\Bbb R^m, \; m\ge 2,$ is the main subject of this article. 

The paper is organized as follows. In Section 2 we introduce the general Markov random flight in the Euclidean spaces $\Bbb R^m, \; m\ge 2$, with arbitrary dissipation function and describe the structure of its distribution. Some basic properties of the joint, conditional and unconditional characteristic functions of the process are also given. In Section 3 we derive the recurrent relations for the joint and conditional densities of the process and of the number of changes of direction in the form of a double convolutions in space and time variables. Basing on these recurrent relations, an integral equation for the transition density of the process is obtained in Section 4, whose solution is given in the form of a uniformly converging series composed of the multiple double convolutions of the singular component of the density with itself. This solution is unique in the class of finitary functions in $\Bbb R^m$. Two important particular cases of the uniform distribution on $S^m(\bold 0,1)$ and of the Gaussian distributions on the unit circumference $S^2(\bold 0,1)$ are considered in Section 5.

\section{Description of Process and Its Basic Properties}

\numberwithin{equation}{section}

Consider the following stochastic motion. A particle starts from the origin $\bold 0 = (0, \dots, 0)$ of the Euclidean space $\Bbb R^m, \; m\ge 2,$ at the initial time instant $t=0$ and moves with some constant speed $c$ (note that $c$ is treated as the constant norm of the velocity). The initial direction is a random $m$-dimensional vector with arbitrary distribution (also called the dissipation function) on the unit sphere 
$$S^m(\bold 0, 1) = \left\{ \bold x=(x_1, \dots ,x_m)\in \Bbb R^m: \; \Vert\bold x\Vert^2 = x_1^2+ \dots +x_m^2=1 \right\} $$
having the absolutely continuous bounded density $\chi(\bold x), \; \bold x\in S^m(\bold 0, 1)$. 
Emphasize that here and thereafter the upper index $m$ means the dimension of the space, which the sphere $S^m(\bold 0, 1)$ is considered in, but not its own dimension which, clearly, is $m-1$.   
The motion is controlled by a homogeneous Poisson process $N(t)$ of rate $\lambda>0$ as follows. At each Poissonian instant, the particle instantaneously takes on a new random direction on $S^m(\bold 0, 1)$ with the same density $\chi(\bold x), \; \bold x\in S^m(\bold 0, 1),$ independently of its previous motion and keeps moving with the same speed $c$ until the next Poisson event occurs, then it takes a new random direction again and so on. 

Let $\bold X(t)=(X_1(t), \dots ,X_m(t))$ be the particle's position at time $t>0$ which is referred to as the $m$-dimensional random flight. At arbitrary time instant $t>0$ the particle, with probability 1, is located in the close $m$-dimensional ball of radius $ct$ 
centered at the origin $\bold 0$:
$$\bold B^m(\bold 0, ct) = \left\{ \bold x=(x_1, \dots ,x_m)\in \Bbb R^m : \; \Vert\bold x\Vert^2 = x_1^2+ \dots +x_m^2\le c^2t^2 \right\} .$$

Consider the probability distribution function 
$$\Phi(\bold x, t) = \text{Pr} \left\{ \bold X(t)\in d\bold x \right\}, \qquad \bold x\in\bold B^m(\bold 0, ct), \quad t\ge 0,$$ 
of process $\bold X(t)$, where $d\bold x$ is the infinitesimal element in the space $\Bbb R^m$ with Lebesgue measure 
$\mu(d\bold x) = dx_1 \dots dx_m$. For arbitrary fixed $t>0$, distribution $\Phi(\bold x, t)$ consists of two components. 

The singular component corresponds to the case when no one Poisson event occurs in the time interval $(0,t)$ and it is concentrated on the sphere
$$S^m(\bold 0, ct) =\partial\bold B^m(\bold 0, ct) = \left\{ \bold x=(x_1, \dots ,x_m)\in \Bbb R^m: \; \Vert\bold x\Vert^2 = x_1^2+ \dots +x_m^2=c^2t^2 \right\} .$$
In this case the particle is located on the sphere $S^m(\bold 0, ct)$ and the probability of this event is
$$\text{Pr} \left\{ \bold X(t)\in S^m(\bold 0, ct) \right\} = e^{-\lambda t} .$$

If at least one Poisson event occurs by time instant $t$, then the particle is located strictly inside the ball $\bold B^m(\bold 0, ct)$ 
and the probability of this event is
$$\text{Pr} \left\{ \bold X(t)\in \text{int} \; \bold B^m(\bold 0, ct) \right\} = 1 - e^{-\lambda t} .$$
The part of the distribution $\Phi(\bold x, t)$ corresponding to this case is concentrated in the interior 
$$\text{int} \; \bold B^m(\bold 0, ct) = \left\{ \bold x=(x_1, \dots ,x_m)\in \Bbb R^m: \; \Vert\bold x\Vert^2 = x_1^2+ \dots +x_m^2<c^2t^2 \right\}$$
of the ball $\bold B^m(\bold 0, ct)$ and forms its absolutely continuous component. 

Let $p(\bold x,t) = p(x_1, \dots ,x_m,t), \;\; \bold x\in\bold B^m(\bold 0, ct) , \; t>0,$ be the density of distribution $\Phi(\bold x, t)$. 
It has the form 
\begin{equation}\label{dens}
p(\bold x,t) = p_s(\bold x,t) + p_{ac}(\bold x,t) , \qquad \bold x\in\bold B^m(\bold 0, ct), \quad t>0,
\end{equation}
where $p_s(\bold x,t)$ is the density (in the sense of generalized functions) of the singular component of $\Phi(\bold x, t)$ concentrated on the sphere $S^m(\bold 0, ct)$ 
and $p_{ac}(\bold x,t)$ is the density of the absolutely continuous component of $\Phi(\bold x, t)$ concentrated in $\text{int} \; \bold B^m(\bold 0, ct)$. 

The density $\chi(\bold x), \;  \bold x\in S^m(\bold 0, 1),$ on the unit sphere $S^m(\bold 0, 1)$ generates the absolutely continuous and bounded (in $\bold x$ for any fixed $t$) density $\varrho(\bold x,t), \bold x\in S^m(\bold 0, ct),$ on the sphere $S^m(\bold 0, ct)$ of radius $ct$ according to the formula 
$\varrho(\bold x,t) = \chi(\frac{1}{ct}\bold x), \; \bold x\in S^m(\bold 0, ct), \; t>0$. Therefore, the singular part of density (\ref{dens}) has the form:
\begin{equation}\label{densS}
p_s(\bold x,t) =  e^{-\lambda t} \varrho(\bold x,t) \delta(c^2t^2-\Vert\bold x\Vert^2) , \qquad t>0,
\end{equation}
where $\delta(x)$ is the Dirac delta-function. 

The absolutely continuous part of density (\ref{dens}) has the form:
\begin{equation}\label{densAC}
p_{ac}(\bold x,t) = f_{ac}(\bold x,t) \Theta(ct-\Vert\bold x\Vert) , \qquad t>0,
\end{equation}
where $f_{ac}(\bold x,t)$ is some positive function absolutely continuous in $\text{int} \; \bold B^m(\bold 0, ct)$ and $\Theta(x)$ is the Heaviside step function given by 
\begin{equation}\label{heaviside}
\Theta(x) = \left\{ \aligned 1, \qquad  & \text{if} \; x>0,\\
                               0, \qquad & \text{if} \; x\le 0.
\endaligned \right.
\end{equation}

Consider the conditional densities $p_n(\bold x,t), \; n\ge 0,$ of process $\bold X(t)$ conditioned by the random events $\{ N(t)=n \}, \; n\ge 0,$ where, remind, $N(t)$  
is the number of the Poisson events that have occurred in the time interval $(0, t)$. Obviously, $p_0(\bold x,t)=\varrho(\bold x,t) \delta(c^2t^2-\Vert\bold x\Vert^2)$. 
Therefore, our aim is to examine conditional densities  $p_n(\bold x,t)$ for $n\ge 1$. 

Consider the conditional characteristic functions (Fourier transform) of process $\bold X(t)$: 
\begin{equation}\label{rec1}
G_n(\boldsymbol\alpha, t) = E \left\{ e^{i\langle\boldsymbol\alpha, \bold X(t)\rangle} \vert \; N(t)=n \right\} , \qquad n\ge 1, 
\end{equation}
where $\boldsymbol\alpha =(\alpha_1, \dots ,\alpha_m) \in \Bbb R^m$ is the real $m$-dimensional vector of inversion parameters and 
$\langle\boldsymbol\alpha, \bold X(t)\rangle$ is the inner product of the vectors $\boldsymbol\alpha$ and $\bold X(t)$.

According to \cite[formula (6.8)]{kol3}, functions (\ref{rec1}) are given by the formula: 
\begin{equation}\label{rec2}
G_n(\boldsymbol\alpha, t) = \frac{n!}{t^n} \int_0^t d\tau_1 \int_{\tau_1}^t d\tau_2 \dots \int_{\tau_{n-1}}^t d\tau_n \left\{ \prod_{j=1}^{n+1}
\psi(\boldsymbol\alpha, \tau_j-\tau_{j-1}) \right\} , \qquad n\ge 1, 
\end{equation}
where  
\begin{equation}\label{rec3}
\psi(\boldsymbol\alpha, t) = \mathcal F_{\bold x} \left[ \varrho(\bold x,t) \delta(c^2t^2-\Vert\bold x\Vert^2) \right](\boldsymbol\alpha) = 
\int_{S^m(\bold 0, ct)} e^{i\langle\boldsymbol\alpha, \bold x\rangle} \; \varrho(\bold x, t) \; \nu(d\bold x)  
\end{equation}
is the characteristic function (Fourier transform) of density $\varrho(\bold x,t)$ concentrated on the  
sphere $S^m(\bold 0, ct)$ of radius $ct$ and $\nu(d\bold x)$ is the surface Lebesgue measure on $S^m(\bold 0, ct)$. 

Consider separately the integral factor in (\ref{rec2}):  
\begin{equation}\label{rec4}
\mathcal J_n(\boldsymbol\alpha, t) = \int_0^t d\tau_1 \int_{\tau_1}^t d\tau_2 \dots
\int_{\tau_{n-1}}^t d\tau_n \left\{ \prod_{j=1}^{n+1} \psi(\boldsymbol\alpha, \tau_j-\tau_{j-1}) \right\} , \qquad n\ge 1. 
\end{equation} 
This function has a quite definite probabilistic sense, namely 
\begin{equation}\label{recc4}
\tilde G_n(\boldsymbol\alpha, t) = \mathcal F_{\bold x} \bigl[ \tilde p_n(\bold x, t) \bigr](\boldsymbol\alpha) = \frac{(\lambda t)^n \; e^{-\lambda t}}{n!} \; G_n(\boldsymbol\alpha, t) = \lambda^n e^{-\lambda t} \mathcal J_n(\boldsymbol\alpha, t), 
\end{equation}
$$\boldsymbol\alpha\in\Bbb R^m, \qquad t>0, \qquad n\ge 1,$$
is the characteristic function (Fourier transform) of the joint probability density $p_n(\bold x, t)$ of the particle's position at time 
instant $t$ and of the number $N(t)=n$ of the Poisson events that have occurred by this time moment $t$. 

It is known (see \cite[Theorem 5]{kol3}) that, for arbitrary $n\ge 1$, functions (\ref{rec4}) are connected with each other 
by the following recurrent relation:  
\begin{equation}\label{rec5}
\mathcal J_n(\boldsymbol\alpha, t) = \int_0^t \psi(\boldsymbol\alpha, t-\tau) \; \mathcal J_{n-1}(\boldsymbol\alpha, \tau) \; d\tau 
= \int_0^t \psi(\boldsymbol\alpha, \tau) \; \mathcal J_{n-1}(\boldsymbol\alpha, t-\tau) \; d\tau, \qquad n\ge 1, 
\end{equation}
where, by definition, ${\mathcal J}_0(\boldsymbol\alpha, t) \overset{\text{def}}{=} \psi(\boldsymbol\alpha, t)$. Formula (\ref{rec5}) can also be represented 
in the following time-convolutional form:
\begin{equation}\label{rec6}
\mathcal J_n(\boldsymbol\alpha, t) = \psi(\boldsymbol\alpha, t) \overset{t}{\ast} \mathcal J_{n-1}(\boldsymbol\alpha, t), \qquad n\ge 1, 
\end{equation} 
where the symbol $\overset{t}{\ast}$ means the convolution operation with respect to time variable $t$. 

From (\ref{rec6}) it follows that 
\begin{equation}\label{rec7}
\mathcal J_n(\boldsymbol\alpha, t) = \left[ \psi(\boldsymbol\alpha, t) \right]^{\overset{t}{\ast}(n+1)} , \qquad n\ge 1,  
\end{equation} 
where $\overset{t}{\ast}(n+1)$ means the $(n+1)$-multiple convolution in $t$.
Applying Laplace transformation $\mathcal L_t$ (in time variable $t$) to (\ref{rec7}), we arrive at the formula 
\begin{equation}\label{rec8}
\mathcal L_t \left[ \mathcal J_n(\boldsymbol\alpha, t) \right] (s) = \bigl( \mathcal L_t \left[ \psi(\boldsymbol\alpha, t) \right] (s) \bigr)^{n+1} , 
\qquad n\ge 1. 
\end{equation} 

It is also known (see \cite[Corollary 5.3]{kol3}) that conditional characteristic functions (\ref{rec2}) satisfy the following recurrent
relation
\begin{equation}\label{rec8}
G_n(\boldsymbol\alpha, t) = \frac{n}{t^n} \int_0^t \tau^{n-1} \psi(\boldsymbol\alpha, t-\tau) \; G_{n-1}(\boldsymbol\alpha, \tau) \; d\tau, \qquad 
G_0(\boldsymbol\alpha, t) \overset{\text{def}}{=} \psi(\boldsymbol\alpha, t), \quad n\ge 1.  
\end{equation} 

The unconditional characteristic function 
\begin{equation}\label{rec9} 
G(\boldsymbol\alpha, t) = E \left\{ e^{i\langle\boldsymbol\alpha, \bold X(t)\rangle} \right\} 
\end{equation} 
of process $\bold X(t)$ satisfies the Volterra integral equation of second kind (see \cite[Theorem 6]{kol3}): 
\begin{equation}\label{rec10}
G(\boldsymbol\alpha, t) = e^{-\lambda t} \psi(\boldsymbol\alpha, t) + \lambda \int_0^t e^{-\lambda (t-\tau)} \psi(\boldsymbol\alpha, t-\tau) \; 
G(\boldsymbol\alpha, \tau) \; d\tau , \qquad t\ge 0, 
\end{equation} 
or in the convolutional form 
\begin{equation}\label{rec11}
G(\boldsymbol\alpha, t) = e^{-\lambda t} \psi(\boldsymbol\alpha, t) + \lambda \bigl[ \left( e^{-\lambda t} \psi(\boldsymbol\alpha, t) \right) 
\overset{t}{\ast} \mathcal J(\boldsymbol\alpha, t) \bigr] .
\end{equation} 
This is the renewal equation for Markov random flight $\bold X(t)$. 

In the class of continuous functions integral equation (\ref{rec10}) (or (\ref{rec11})) has the unique solution given by the uniformly converging series 
\begin{equation}\label{rec12} 
G(\boldsymbol\alpha, t) = e^{-\lambda t} \sum_{n=0}^{\infty} \lambda^n \; \left[ \psi(\boldsymbol\alpha, t) \right]^{\overset{t}{\ast} (n+1)} \; .
\end{equation} 

From (\ref{rec11}) we obtain the general formula for the Laplace transform of characteristic function (\ref{rec9}): 
\begin{equation}\label{rec13}
\mathcal L_t \left[ G(\boldsymbol\alpha, t) \right] (s) = \frac{\mathcal L_t \left[ \psi(\boldsymbol\alpha, t) \right] (s+\lambda)}{1 - \lambda \; 
\mathcal L_t \left[ \psi(\boldsymbol\alpha, t) \right](s+\lambda)} \; , \qquad \text{Re} \; s > 0.
\end{equation} 

These properties will be used in the next section for deriving recurrent relations for the joint and conditional densities of Markov 
random flight $\bold X(t)$.

\section{Recurrent Relations}

\numberwithin{equation}{section}

Consider the joint probability densities $p_n(\bold x, t), \; n\ge 0, \; \bold x\in\bold B^m(\bold 0, ct), \; t>0,$ 
of the particle's position $\bold X(t)$ at time instant $t>0$ and of the number of the Poisson events $\{ N(t)=n\}$ that have occurred by this instant $t$. 
For $n=0$, we have 
\begin{equation}\label{joint0}
p_0(\bold x, t) = p_s(\bold x, t) = e^{-\lambda t} \varrho(\bold x,t) \delta(c^2t^2-\Vert\bold x\Vert^2) , \qquad t>0,  
\end{equation} 
where, remind, $p_s(\bold x, t)$ is the singular part of density (\ref{dens}) concentrated on the surface of the sphere 
$S^m(\bold 0, ct) = \partial\bold B^m(\bold 0, ct)$ and given by (\ref{densS}). 

If $n\ge 1$, then, according to (\ref{densAC}), joint densities $p_n(\bold x, t)$ 
have the form: 
\begin{equation}\label{jointN}
p_n(\bold x, t) = f_n(\bold x,t) \Theta(ct-\Vert\bold x\Vert) , \qquad n\ge 1, \quad t>0, 
\end{equation} 
where $f_n(\bold x,t), \; n\ge 1,$ are some positive functions absolutely continuous in $\text{int} \; \bold B^m(\bold 0, ct)$ and $\Theta(x)$ is 
the Heaviside step function. 

The joint density $p_{n+1}(\bold x,t)$ can be expressed through the previous one $p_n(\bold x,t)$ 
by means of a recurrent relation. This result is given by the following theorem.

\bigskip

{\bf Theorem 1.} {\it The joint densities} $p_n(\bold x, t), \; n\ge 1,$ {\it are connected with 
each other by the following recurrent relation:} 
\begin{equation}\label{rec14}
p_{n+1}(\bold x, t) = \lambda \int_0^t \bigl[ p_0(\bold x, t-\tau) \overset{\bold x}{\ast} p_n(\bold x, \tau) \bigr] \; d\tau , 
\qquad n\ge 1, \quad \bold x\in\text{int} \; \bold B^m(\bold 0, ct), \quad t>0.
\end{equation} 

\vskip 0.2cm 

\begin{proof}
Applying Fourier transformation to the right-hand side of (\ref{rec15}), we have: 

\begin{equation}\label{rec17}
\aligned 
\lambda \; \mathcal F_{\bold x} & \left[ \int_0^t \bigl[ p_0(\bold x, t-\tau) \overset{\bold x}{\ast} p_n(\bold x, \tau) \bigr] \; d\tau   \right](\boldsymbol\alpha) \\ 
& = \lambda \int_0^t \mathcal F_{\bold x} \left[ p_0(\bold x, t-\tau) \overset{\bold x}{\ast} p_n(\bold x, \tau) \right](\boldsymbol\alpha) \; d\tau \\
& = \lambda \int_0^t \mathcal F_{\bold x} \bigl[ p_0(\bold x, t-\tau) \bigr](\boldsymbol\alpha) \; \mathcal F_{\bold x} \bigl[ p_n(\bold x, \tau) \bigr](\boldsymbol\alpha) \; d\tau \\
& = \lambda \int_0^t e^{-\lambda(t-\tau)} \mathcal F_{\bold x} \bigl[ \varrho(\bold x, t-\tau) \delta(c^2(t-\tau)^2-\Vert\bold x\Vert^2) \bigr](\boldsymbol\alpha) \;\; \mathcal F_{\bold x} \bigl[ p_n(\bold x, \tau) \bigr](\boldsymbol\alpha) \; d\tau \\
& = \lambda \int_0^t e^{-\lambda(t-\tau)} \psi(\boldsymbol\alpha, t-\tau) \; \lambda^n e^{-\lambda\tau} \mathcal J_n(\boldsymbol\alpha, \tau) \; d\tau \\ 
& = \lambda^{n+1} e^{-\lambda t} \int_0^t \psi(\boldsymbol\alpha, t-\tau) \; \mathcal J_n(\boldsymbol\alpha, \tau) \; d\tau \\ 
& = \lambda^{n+1} e^{-\lambda t} \mathcal J_{n+1}(\boldsymbol\alpha, t) \\
& = \mathcal F_{\bold x} \bigl[ p_{n+1}(\bold x, t) \bigr](\boldsymbol\alpha) , 
\endaligned
\end{equation}
where we have used formulas (\ref{rec3}), (\ref{recc4}), (\ref{rec5}). Thus, both the functions on the left- and right-hand sides of (\ref{rec14}) have the same Fourier transform and, therefore, they coincide. 

The change of integration order in the first step of (\ref{rec17}) is justified because the convolution $p_0(\bold x, t-\tau) \overset{\bold x}{\ast} p_n(\bold x, \tau)$ of the singular part $p_0(\bold x, t-\tau)$ of the density with the absolutely continuous one  $p_n(\bold x, \tau), \; n\ge 1,$ is the absolutely continuous (and, therefore, uniformly bounded in $\bold x$) function. From this fact it follows that, for any $n\ge 1$, the integral in square brackets on the left-hand side of (\ref{rec17}) converges uniformly in $\bold x$ for any fixed $t$. The theorem is proved.  
\end{proof}

\bigskip

{\it Remark 1.} In view of (\ref{densS}) and (\ref{densAC}), formula (\ref{rec14}) can be represented in the following expanded form: 
\begin{equation}\label{rec15}
\aligned 
p_{n+1}(\bold x, t) & = \lambda \int_0^t e^{-\lambda(t-\tau)} \\
& \quad \times \left\{ \int \varrho(\bold x-\boldsymbol\xi,t-\tau) \delta(c^2(t-\tau)^2-\Vert\bold x-\boldsymbol\xi\Vert^2) 
\; f_n(\boldsymbol\xi, \tau) \Theta(c\tau-\Vert\boldsymbol\xi\Vert) \; \nu(d\boldsymbol\xi) \right\} d\tau ,
\endaligned
\end{equation}
$$n\ge 1, \quad \bold x\in\text{int} \; \bold B^m(\bold 0, ct), \quad t>0,$$
where the function $f_n(\boldsymbol\xi, \tau)$ is absolutely continuous in the variable $\boldsymbol\xi=(\xi_1,\dots,\xi_m)\in\Bbb R^m$ and $\nu(d\boldsymbol\xi)$ 
is the surface Lebesgue measure. Integration area in the interior integral on the right-hand side of (\ref{rec15}) is determined by all the $\boldsymbol\xi$, 
under which the integrand takes non-zero values, that is, by the system 
$$\boldsymbol\xi\in\Bbb R^m  \; : \; \left\{ \aligned & \Vert\bold x-\boldsymbol\xi\Vert^2 = c^2(t-\tau)^2 , \\
                                                      & \Vert\boldsymbol\xi\Vert < c\tau. \endaligned \right.$$
The first relation of this system determines a sphere $S^m(\bold x, c(t-\tau))$ of radius $c(t-\tau)$ centred at point $\bold x$, while the second one represents 
an open ball $\text{int} \; \bold B^m(\bold 0, c\tau)$ of radius $c\tau$ centred at the origin $\bold 0$. Their intersection 
\begin{equation}\label{setM}
M(\bold x, \tau) = S^m(\bold x, c(t-\tau))\cap\text{int} \; \bold B^m(\bold 0, c\tau),
\end{equation} 
which is a part of (or the whole) surface of sphere $S^m(\bold x, c(t-\tau))$ located inside the ball $\bold B^m(\bold 0, c\tau)$, represents the integration area of dimension $m-1$ in the interior integral of (\ref{rec15}). Note that the sum of the radii of $S^m(\bold x, c(t-\tau))$ and $\text{int} \; \bold B^m(\bold 0, c\tau)$ is $c(t-\tau)+c\tau = ct > \Vert\bold x\Vert$, that is greater than the distance $\Vert\bold x\Vert$ between their centres $\bold 0$ and $\bold x$. This fact, as well as some simple geometric reasonongs, show that intersection (\ref{setM}) depends on $\tau\in (0,t)$ as follows. 

$\bullet$ If $\tau\in (0, \; \frac{t}{2} - \frac{\Vert x\Vert}{2c}]$, then intersection (\ref{setM}) is empty, that is, 
$M(\bold x, \tau) = \varnothing$.

$\bullet$ If $\tau\in (\frac{t}{2} - \frac{\Vert x\Vert}{2c}, \; \frac{t}{2} + \frac{\Vert x\Vert}{2c}]$, then intersection 
$M(\bold x, \tau)$ is not empty and represents some hyperspace of dimension $m-1$. 

$\bullet$ If $\tau\in (\frac{t}{2} + \frac{\Vert x\Vert}{2c}, \; t]$, then $S^m(\bold x, c(t-\tau))\subset \text{int} \; \bold B^m(\bold 0, c\tau)$ and, therefore, in this case $M(\bold x, \tau) = S^m(\bold x, c(t-\tau))$. 

Thus, formula (\ref{rec15}), as well as (\ref{rec14}), can be rewritten in the expanded form: 
\begin{equation}\label{rec16}
\aligned 
p_{n+1}(\bold x, t) & = \lambda \int\limits_{\frac{t}{2} - \frac{\Vert x\Vert}{2c}}^{\frac{t}{2} + \frac{\Vert x\Vert}{2c}} e^{-\lambda(t-\tau)} \left\{ \int\limits_{M(\bold x, \tau)} \varrho(\bold x-\boldsymbol\xi,t-\tau) \; f_n(\boldsymbol\xi, \tau) \; \nu(d\boldsymbol\xi) \right\} d\tau \\ 
& + \lambda \int\limits_{\frac{t}{2} + \frac{\Vert x\Vert}{2c}}^t e^{-\lambda(t-\tau)} \left\{ \int\limits_{S^m(\bold x, c(t-\tau))} \varrho(\bold x-\boldsymbol\xi,t-\tau) \; f_n(\boldsymbol\xi, \tau) \; \nu(d\boldsymbol\xi) \right\} d\tau
\endaligned 
\end{equation}
and the expressions in curly brackets of (\ref{rec16}) represents surface integrals over $M(\bold x, \tau)$ and $S^m(\bold x, c(t-\tau))$. 

\bigskip 

{\it Remark 2.} By means of the double convolution operation of two arbitrary generalized functions $f_1(\bold x, t), f_2(\bold x, t) \in\mathscr{S'}, \; \bold x\in\Bbb R^m, \; t>0,$ 
\begin{equation}\label{rec16a}
f_1(\bold x, t) \overset{\bold x}{\ast} \overset{t}{\ast} f_2(\bold x, t) = \int_0^t \int_{\Bbb R^m} f_1(\boldsymbol\xi, \tau) \; f_2(\bold x - \boldsymbol\xi, t-\tau) \; d\boldsymbol\xi \; d\tau 
\end{equation} 
formula (\ref{rec14}) can be represented in the succinct convolutional form 
\begin{equation}\label{rec16b}
p_{n+1}(\bold x, t) = \lambda \left[ p_0(\bold x, t) \overset{\bold x}{\ast} \overset{t}{\ast} p_n(\bold x, t) \right] . 
\end{equation}

\bigskip 

Taking into account the well-known connections between the joint and conditional densities, we can extract from Theorem 1 a convolution-type recurrent relation for the 
conditional probability densities $\tilde p_n(\bold x, t), \; n\ge 1$. 

\bigskip 

{\bf Corollary 1.1.} {\it The conditional densities} $\tilde p_n(\bold x, t), \; n\ge 1,$ {\it are connected with each other by the following recurrent relation:}  
\begin{equation}\label{rec18}
\tilde p_{n+1}(\bold x, t) = \frac{n+1}{t^{n+1}} \int_0^t \tau^n \bigl[ \tilde p_0(\bold x, t-\tau) \overset{\bold x}{\ast} \; \tilde p_n(\bold x, \tau) \bigr] \; d\tau , 
\qquad n\ge 1, \quad \bold x\in\text{int} \; \bold B^m(\bold 0, ct), \quad t>0, 
\end{equation} 
{\it where $\tilde p_0(\bold x, t) = \varrho(\bold x,t) \delta(c^2t^2-\Vert\bold x\Vert^2)$ is the conditional density corresponding to the 
case when no one Poisson event occurs before time instant $t$.}

\vskip 0.2cm

\begin{proof}
The proof immediately follows from Theorem 1 and recurrent formula (\ref{rec8}).
\end{proof}

\bigskip

{\it Remark 3.} Formulas (\ref{rec14}) and (\ref{rec18}) are also valid for $n=0$. In this case, for arbitrary $t>0$, they take the form: 
\begin{equation}\label{rec19}
p_1(\bold x, t) = \lambda \int_0^t \bigl[ p_0(\bold x, t-\tau) \overset{\bold x}{\ast} \; p_0(\bold x, \tau) \bigr] \; d\tau , 
\end{equation}
\begin{equation}\label{rec20}
\tilde p_1(\bold x, t) = \frac{1}{t} \int_0^t \bigl[ \tilde p_0(\bold x, t-\tau) \overset{\bold x}{\ast} \; \tilde p_0(\bold x, \tau) \bigr] \; d\tau , 
\end{equation}
where, remind, function $p_0(\bold x, t)$ defined by (\ref{joint0}) is the singular part of the density concentrated on the surface 
of the sphere $S^m(\bold 0, ct)$. The derivation of (\ref{rec19}) is a simple recompilation of the proof of Theorem 1 and taking into account 
the boundedness of $p_0(\bold x, t)$ which justifies the change of integration order in (\ref{rec17}). Formula (\ref{rec20}) follows from (\ref{rec19}).

\section{Integral Equation for Transition Density}

\numberwithin{equation}{section}

The transition probability density $p(\bold x, t)$ of the multidimensional Markov flight $\bold X(t)$ is defined by the formula 
\begin{equation}\label{int1}
p(\bold x, t) = \sum_{n=0}^{\infty} p_n(\bold x, t) , \qquad \bold x\in\bold B^m(\bold 0, ct), \quad t>0, 
\end{equation}
where the joint densities $p_n(\bold x, t), \; n\ge 0,$ are given by (\ref{joint0}) and (\ref{jointN}). The density (\ref{int1}) 
is defined everywhere in the ball $\bold B^m(\bold 0, ct)$, while the function 
\begin{equation}\label{int2}
p_{ac}(\bold x, t) = \sum_{n=1}^{\infty} p_n(\bold x, t) 
\end{equation} 
forms its absolutely continuous part concentrated in the interior $\text{int} \; \bold B^m(\bold 0, ct)$ of the ball. 
Therefore, series (\ref{int2}) converges uniformly everywhere in the close ball $\bold B^m(\bold 0, ct-\varepsilon)$ 
for arbitrary small $\varepsilon>0$. 

In the following theorem we state an integral equation for density (\ref{int1}). 

\bigskip

{\bf Theorem 2.} {\it The transition probability density $p(\bold x, t)$ of the Markov random flight $\bold X(t)$ 
satisfies the integral equation:} 
\begin{equation}\label{int3} 
p(\bold x, t) = p_0(\bold x, t) + \lambda \int_0^t \bigl[ p_0(\bold x, t-\tau) \overset{\bold x}{\ast} \; p(\bold x, \tau) \bigr] \; d\tau , \qquad \bold x\in\bold B^m(\bold 0, ct), \quad t>0.
\end{equation} 

{\it In the class of finitary functions (that is, generalized functions defined on the compact sets of $\Bbb R^m$), integral equation} (\ref{int3}) 
{\it has the unique solution given by the series}  
\begin{equation}\label{int4}
p(\bold x,t) = \sum_{n=0}^{\infty} \lambda^n \left[ p_0(\bold x, t) \right]^{\overset{\bold x}{\ast} \overset{t}{\ast}(n+1)} , 
\end{equation}
{\it where the symbol $\overset{\bold x}{\ast} \overset{t}{\ast}(n+1)$ means the $(n+1)$-multiple double convolution with respect to spatial and time variables defined by} (\ref{rec16a}), {\it that is,}
$$\left[ p_0(\bold x, t) \right]^{\overset{\bold x}{\ast} \overset{t}{\ast}(n+1)} = \underbrace{p_0(\bold x, t) \overset{\bold x}{\ast} \overset{t}{\ast} p_0(\bold x, t) \overset{\bold x}{\ast} \overset{t}{\ast} \dots \overset{\bold x}{\ast} \overset{t}{\ast} p_0(\bold x, t)}_{(n+1) \; \text{times}} .$$
{\it Series} (\ref{int4}) {\it is convergent everywhere in the open ball} $\text{int} \; \bold B^m(\bold 0, ct)$.  
{\it For any small $\varepsilon>0$, the series} (\ref{int4}) {\it converges uniformly (in $\bold x$ for any fixed $t>0$) in the close ball $\bold B^m(\bold 0, ct-\varepsilon)$ and, therefore, it determines the density $p(\bold x, t)$ which is continuous and bounded in this ball}.  

\vskip 0.2cm 

\begin{proof}
Applying Theorem 1 and taking into account the uniform convergence of series (\ref{int2}) and of the integral in formula (\ref{rec14}), we have: 

$$\aligned 
p(\bold x, t) & = \sum_{n=0}^{\infty} p_n(\bold x, t) \\
& = p_0(\bold x, t) + \sum_{n=1}^{\infty} p_n(\bold x, t) \\ 
& = p_0(\bold x, t) + \lambda \sum_{n=1}^{\infty} \int_0^t \bigl[ p_0(\bold x, t-\tau) \overset{\bold x}{\ast} \; p_{n-1}(\bold x, \tau) \bigr] \; d\tau \\ 
& = p_0(\bold x, t) + \lambda \int_0^t \sum_{n=1}^{\infty} \bigl[ p_0(\bold x, t-\tau) \overset{\bold x}{\ast} \; p_{n-1}(\bold x, \tau) \bigr] \; d\tau \\ 
& = p_0(\bold x, t) + \lambda \int_0^t \left[ p_0(\bold x, t-\tau) \overset{\bold x}{\ast} \; \left\{ \sum_{n=1}^{\infty} p_{n-1}(\bold x, \tau) \right\} \right] \; d\tau \\ 
& = p_0(\bold x, t) + \lambda \int_0^t \left[ p_0(\bold x, t-\tau) \overset{\bold x}{\ast} \; \left\{ \sum_{n=0}^{\infty} p_n(\bold x, \tau) \right\} \right] \; d\tau \\ 
& = p_0(\bold x, t) + \lambda \int_0^t \bigl[ p_0(\bold x, t-\tau) \overset{\bold x}{\ast} \; p(\bold x, \tau) \bigr] \; d\tau , 
\endaligned$$
proving (\ref{int3}). 

Another way of proving the theorem is to apply the Fourier transformation to both sides of (\ref{int3}). Justifying then the change of the order of integrals similarly as it was done in (\ref{rec17}), we arrive at Volterra integral equation (\ref{rec10}) for Fourier transforms. 

Using notation (\ref{rec16a}), equation (\ref{int3}) can be represented in the convolutional form 
\begin{equation}\label{int5}
p(\bold x, t) = p_0(\bold x, t) + \lambda \bigl[ p_0(\bold x, t) \overset{\bold x}{\ast} \overset{t}{\ast} p(\bold x, t) \bigr] , 
\qquad \bold x\in\bold B^m(\bold 0, ct), \quad t>0.
\end{equation}
Let us check that series (\ref{int4}) satisfies equation (\ref{int5}). Substituting (\ref{int4}) into the right-hand side of (\ref{int5}), we have: 

$$\aligned 
p_0(\bold x, t) + \lambda \biggl[ p_0(\bold x, t) \overset{\bold x}{\ast} \overset{t}{\ast} \biggl( \sum_{n=0}^{\infty} \lambda^n \left[ p_0(\bold x, t) \right]^{\overset{\bold x}{\ast} \overset{t}{\ast}(n+1)} \biggr) \biggr] & = p_0(\bold x, t) + \sum_{n=0}^{\infty} \lambda^{n+1} \left[ p_0(\bold x, t) \right]^{\overset{\bold x}{\ast} \overset{t}{\ast}(n+2)} \\ 
& = p_0(\bold x, t) + \sum_{n=1}^{\infty} \lambda^n \left[ p_0(\bold x, t) \right]^{\overset{\bold x}{\ast} \overset{t}{\ast}(n+1)} \\ 
& = \sum_{n=0}^{\infty} \lambda^n \left[ p_0(\bold x, t) \right]^{\overset{\bold x}{\ast} \overset{t}{\ast}(n+1)} \\ 
& = p(\bold x, t) 
\endaligned$$
and, therefore, series (\ref{int4}) is really the solution to equation (\ref{int5}). 

Note that applying Fourier transformation to (\ref{int3}) and (\ref{int4}) and taking into account (\ref{densS}), we arrive at the known results (\ref{rec11}) and (\ref{rec12}), respectively. The uniqueness of solution (\ref{int4}) in the class of finitary functions follows from the uniqueness of the solution of Volterra integral equation (\ref{rec10}) for its Fourier transform (\ref{rec12}) (i.e. characteristic function) in the class of continuous functions. 

Since the transition density $p(\bold x,t)$ is absolutely continuous in the open ball $\text{int} \; \bold B^m(\bold 0, ct)$, then, for any $\varepsilon>0$, it is continuous and uniformly bounded in the close ball $\bold B^m(\bold 0, ct-\varepsilon)$. From this fact and taking into account the uniqueness of the solution of integral equation (\ref{int3}) in the class of finitary functions, we can conclude that series (\ref{int4}) converges uniformly in $\bold B^m(\bold 0, ct-\varepsilon)$ for any small $\varepsilon>0$. This completes the proof. 
\end{proof}

\section{Some Particular Cases}

\numberwithin{equation}{section}

In this section we consider two important particular cases of the general Markov random flight described in Section 2 when the dissipation function has the uniform distribution on the unit sphere $S^m(\bold 0,1)$ and Gaussian distribution on the unit circumference $S^2(\bold 0,1)$.

\subsection{Symmetric Random Flights} 

Suppose that the initial and every new direction are chosen according to the uniform distribution on the unit sphere $S^m(\bold 0,1)$. Such processes in the Euclidean spaces $\Bbb R^m$ of different dimensions $m\ge 2$, which are referred to as the symmetric Markov random flights, have become the subject of a series of works \cite{kol2,kol3,kol4,kol5,kol6,kol7}, \cite{mas}, \cite{sta1,sta2}. 

In this symmetric case the function $\varrho(\bold x,t)$ is the density of the uniform distribution on the surface of the sphere 
$S^m(\bold 0, ct)$ and, therefore, it does not depend on spatial variable $\bold x$. Then, according to (\ref{densS}), the singular part of the transition density of process $\bold X(t)$ takes the form: 
\begin{equation}\label{sym1}
p_s(\bold x,t) =  e^{-\lambda t} \frac{\Gamma\left( \frac{m}{2} \right)}{2\pi^{m/2}\; (ct)^{m-1}} \; \delta(c^2t^2-\Vert\bold x\Vert^2) , \qquad m\ge 2, \quad t>0.
\end{equation}

Therefore, according to Theorem 1, for arbitrary dimension $m\ge 2$, the absolutely continuous parts $f_n(\bold x, t)$, \; $n\ge 0,$  of the joint probability densities of the symmetric Markov random flight are connected with each other by the following recurrent relation: 
\begin{equation}\label{sym2}
f_{n+1}(\bold x, t) = \frac{\lambda \; \Gamma\left( \frac{m}{2} \right)}{2\pi^{m/2}\; c^{m-1}}  \int_0^t \frac{e^{-\lambda(t-\tau)}}{(t-\tau)^{m-1}} \biggl\{ \int\limits_{M(\bold x,\tau)} f_n(\boldsymbol\xi,\tau) \; d\boldsymbol\xi \biggr\} \; d\tau , 
\end{equation}
$$\bold x=(x_1,\dots,x_m) \in\text{int} \; \bold B^m(\bold 0, ct), \quad m\ge 2, \quad n\ge 0, \quad t>0,$$
where the integration area $M(\bold x,\tau)$ is given by (\ref{setM}).

It is known (see \cite[formula (7)]{kol4}) that, in arbitrary dimension $m\ge 2$, the joint density of symmetric Markov random flight $\bold X(t)$ and of the single change of direction is given by the formula 
\begin{equation}\label{sym3}
f_1(\bold x,t) = \lambda e^{-\lambda t} \; \frac{2^{m-3} \Gamma\left( \frac{m}{2} \right)}{\pi^{m/2} c^m t^{m-1}} \; F\left( \frac{m-1}{2}, -\frac{m}{2}+2; \; \frac{m}{2}; \; \frac{\Vert\bold x\Vert^2}{c^2t^2} \right) , 
\end{equation}
$$\bold x=(x_1,\dots,x_m)\in\text{int} \; \bold B^m(\bold 0, ct), \quad m\ge 2, \quad t>0,$$
where 
$$F(\alpha,\beta;\gamma;z) = \sum_{k=0}^{\infty} \frac{(\alpha)_k (\beta)_k}{(\gamma)_k} \; \frac{z^k}{k!}$$
is the Gauss hypergeometric function. 

Then, by substituting (\ref{sym3}) into (\ref{sym2}) (for $n=1$), we obtain the following formula for the joint density of process $\bold X(t)$ and of two changes of direction:
\begin{equation}\label{sym4}
\aligned
f_2(\bold x,t) & = \lambda^2 e^{-\lambda t} \; \frac{2^{m-4} \left[\Gamma\left( \frac{m}{2} \right)\right]^2}{\pi^m \; c^{2m-1}} \\  
& \qquad\times \int_0^t \biggl\{ \int\limits_{M(\bold x,\tau)} F\left( \frac{m-1}{2}, -\frac{m}{2}+2; \; \frac{m}{2}; \; \frac{\Vert\boldsymbol\xi\Vert^2}{c^2\tau^2} \right) d\boldsymbol\xi \biggr\} \; \frac{d\tau}{(\tau(t-\tau))^{m-1}} , 
\endaligned 
\end{equation}
$$\bold x=(x_1,\dots,x_m)\in\text{int} \; \bold B^m(\bold 0, ct), \quad m\ge 2, \quad t>0.$$

In the three-dimensional Euclidean space $\Bbb R^3$, joint density (\ref{sym3}) was computed explicitly by different methods and it has the form (see \cite[formula (25)]{kol5} or \cite[the second term of formulas (1.3) and (4.21)]{sta1}):
\begin{equation}\label{sym5}
f_1(\bold x,t) = \frac{\lambda e^{-\lambda t}}{4\pi c^2 t \Vert\bold x\Vert} \; \ln\left( \frac{ct+\Vert\bold x\Vert}{ct-\Vert\bold x\Vert} \right) , 
\end{equation}
$$\bold x=(x_1,x_2,x_3)\in\text{int} \; \bold B^3(\bold 0, ct), \quad \Vert\bold x\Vert=\sqrt{x_1^2+x_2^2+x_3^2}, \quad t>0.$$
By substituting this joint density into (\ref{sym2}) (for $n=1, \; m=3$), we arrive at the formula:  
\begin{equation}\label{sym6}
f_2(\bold x,t) = \frac{\lambda^2 e^{-\lambda t}}{16 \pi^2 c^4}    
\int_0^t \biggl\{ \int\limits_{M(\bold x,\tau)} \ln\left( \frac{c\tau+\Vert\boldsymbol\xi\Vert}{c\tau-\Vert\boldsymbol\xi\Vert} \right) \frac{d\boldsymbol\xi}{\Vert\boldsymbol\xi\Vert} \biggr\} \; \frac{d\tau}{\tau(t-\tau)^2} , 
\end{equation}
$$\bold x=(x_1,x_2,x_3)\in\text{int} \; \bold B^3(\bold 0, ct), \quad t>0.$$
Formula (\ref{sym6}) can also be obtained by setting $m=3$ in (\ref{sym4}).

According to Theorem 2 and (\ref{sym1}), the transition density of the $m$-dimensional symmetric Markov random flight solves the integral equation 
\begin{equation}\label{sym7} 
\aligned 
p(\bold x, t) & = \frac{\Gamma\left( \frac{m}{2} \right)}{2\pi^{m/2}\; c^{m-1}} \biggl\{ 
\frac{e^{-\lambda t}}{t^{m-1}} \; \delta(c^2t^2-\Vert\bold x\Vert^2) \\ 
& \qquad + \lambda \int_0^t \biggl[ \left( \frac{e^{-\lambda (t-\tau)}}{(t-\tau)^{m-1}} \; \delta(c^2(t-\tau)^2-\Vert\bold x\Vert^2) \right) \overset{\bold x}{\ast} \; p(\bold x, \tau) \biggr] \; d\tau \biggr\} , 
\endaligned
\end{equation}
$$\bold x=(x_1,\dots,x_m)\in\bold B^m(\bold 0, ct), \quad t>0.$$
In the class of finitary functions, equation (\ref{sym7}) has the unique solution given by the series 
\begin{equation}\label{sym8}
p(\bold x,t) = \sum_{n=0}^{\infty} \lambda^n \left( \frac{\Gamma\left( \frac{m}{2} \right)}{2\pi^{m/2} \; c^{m-1}} \right)^{n+1} 
\left[ \frac{e^{-\lambda t}}{t^{m-1}} \; \delta(c^2t^2-\Vert\bold x\Vert^2) \right]^{\overset{\bold x}{\ast} \overset{t}{\ast}(n+1)} . 
\end{equation}

\bigskip

\subsection{Gaussian Distribution on Circumference} 

Consider now the case of the non-symmetric planar random flight when the initial and each new direction are chosen according to the Gaussian distribution on the unit circumference $S^2(\bold 0,1)$ with the two-dimensional density
\begin{equation}\label{gauss1}
\chi_k(\bold x) = \frac{1}{2\pi \; I_0(k)} \; \exp\left( \frac{k x_1}{\Vert\bold x\Vert} \right) \; \delta(1-\Vert\bold x\Vert^2) , 
\end{equation}
$$\bold x =(x_1,x_2)\in\Bbb R^2, \qquad \Vert\bold x\Vert=\sqrt{x_1^2+x_2^2} \qquad k\in\Bbb R,$$
where $I_0(z)$ is the modified Bessel function of order 0. Formula (\ref{gauss1}) determines the one-parametric family of Gaussian densities $\left\{ \chi_k(\bold x), \; k\in\Bbb R \right\}$, and for any fixed real $k\in\Bbb R$ the density $\chi_k(\bold x)$ 
is absolutely continuous and uniformly bounded on $S^2(\bold 0,1)$. If $k=0$, then formula (\ref{gauss1}) transforms into the density of the uniform distribution on the unit circumference $S^2(\bold 0,1)$, while for $k\neq 0$ it produces pure Gaussian-type densities.  

In the unit polar coordinates $x_1=\cos\theta, \; x_2=\sin\theta,$ formula (\ref{gauss1}) takes the form of the circular Gaussian law: 
\begin{equation}\label{gauss2}
\chi_k(\theta) = \frac{e^{k \cos\theta}}{2\pi \; I_0(k)} , \qquad \theta\in [-\pi, \pi), \quad k\in\Bbb R.   
\end{equation}

For arbitrary real $k\in\Bbb R$, Gaussian density (\ref{gauss1}) on the unit circumference $S^2(\bold 0,1)$ generates 
the Gaussian density 
\begin{equation}\label{gauss3}
p_s(\bold x,t) =  \frac{e^{-\lambda t}}{2\pi ct \; I_0(k)} \; \exp\left( \frac{k x_1}{\Vert\bold x\Vert} \right) \;
\delta(c^2t^2-\Vert\bold x\Vert^2) , 
\end{equation}
$$\bold x =(x_1,x_2)\in\Bbb R^2, \quad \Vert\bold x\Vert=\sqrt{x_1^2+x_2^2}, \quad  t>0, \quad k\in\Bbb R,$$
concentrated on the circumference $S^2(\bold 0, ct)$ of radius $ct$.
Then, according to Theorem 1, the joint densities are connected with each other by the recurrent relation 
\begin{equation}\label{gauss4}
\aligned 
& f_{n+1}(\bold x, t) \\  
& = \frac{\lambda}{2\pi c I_0(k)} \int_0^t \biggl\{ \int\limits_{M(\bold x,\tau)} \exp\left( 
\frac{k (x_1-\xi_1)}{\sqrt{(x_1-\xi_1)^2+(x_2-\xi_2)^2}} \right) \; f_n(\xi_1,\xi_2,\tau) \; d\xi_1 d\xi_2 \biggr\} \frac{e^{-\lambda(t-\tau)}}{t-\tau} \; d\tau ,
\endaligned
\end{equation}
$$\bold x =(x_1,x_2)\in \; \text{int} \; \bold B^2(\bold 0, ct), \quad n\ge 0, \quad  t>0, \quad k\in\Bbb R.$$

According to Theorem 2 and (\ref{gauss3}), the transition density of the planar Markov random flight with Gaussian dissipation function (\ref{gauss1}) satisfies the integral equation 
\begin{equation}\label{gauss5}
\aligned 
p(\bold x,t) & = \frac{e^{-\lambda t}}{2\pi ct \; I_0(k)} \; \exp\left( \frac{k x_1}{\Vert\bold x\Vert} \right) \;
\delta(c^2t^2-\Vert\bold x\Vert^2) \\ 
& \qquad + \frac{\lambda}{2\pi c \; I_0(k)} \int_0^t \left[ \left( \frac{e^{-\lambda\tau}}{\tau} \; \exp\left( \frac{k x_1}{\Vert\bold x\Vert} \right) \; \delta(c^2\tau^2-\Vert\bold x\Vert^2) \right) \overset{\bold x}{\ast} \; p(\bold x, \tau) \right] d\tau , 
\endaligned
\end{equation}
$$\bold x=(x_1,x_2)\in\bold B^2(\bold 0, ct), \quad \Vert\bold x\Vert=\sqrt{x_1^2+x_2^2}, \quad t>0, \quad k\in\Bbb R.$$ 
In the class of finitary functions, equation (\ref{gauss5}) has the unique solution given by the series 
\begin{equation}\label{gauss6}
p(\bold x,t) = \sum_{n=0}^{\infty} \lambda^n \left( \frac{1}{2\pi c \; I_0(k)} \right)^{n+1} \left[ \frac{e^{-\lambda t}}{t} \; \exp\left( \frac{k x_1}{\Vert\bold x\Vert} \right) \; \delta(c^2 t^2-\Vert\bold x\Vert^2) \right]^{\overset{\bold x}{\ast} \overset{t}{\ast}(n+1)} . 
\end{equation}

\end{document}